\documentclass[12pt]{article}

\catcode`\@=11

\thicklines
\newskip\Einheit \Einheit=.6cm
\newcount\xcoord \newcount\ycoord
\newdimen\xdim \newdimen\ydim \newdimen\PfadD@cke \newdimen\Pfadd@cke
\def\PfadDicke#1{\PfadD@cke#1 \divide\PfadD@cke by2 
\Pfadd@cke\PfadD@cke \multiply\PfadD@cke by2}
\long\def\LOOP#1\REPEAT{\def\BODY{#1}\ITERATE}
\def\ITERATE{\BODY \let\next\ITERATE \else\let\next\relax\fi \next}
\let\REPEAT=\fi
\def\Punkt{\hbox{\raise-2pt\hbox to0pt{\hss\scriptsize$\bullet$\hss}}}

\def\DuennPunkt(#1,#2){\unskip
  \raise#2 \Einheit\hbox to0pt{\hskip#1 \Einheit
          \raise-1.5pt\hbox to0pt{\hss\tiny$\bullet$\hss}\hss}}
		  
\def\NormalPunkt(#1,#2){\unskip
  \raise#2 \Einheit\hbox to0pt{\hskip#1 \Einheit
          \raise-3pt\hbox to0pt{\hss\large$\bullet$\hss}\hss}}
\def\DickPunkt(#1,#2){\unskip
  \raise#2 \Einheit\hbox to0pt{\hskip#1 \Einheit
          \raise-4pt\hbox to0pt{\hss\Large$\bullet$\hss}\hss}}
\def\Kreis(#1,#2){\unskip
  \raise#2 \Einheit\hbox to0pt{\hskip#1 \Einheit
          \raise-4pt\hbox to0pt{\hss\Large$\circ$\hss}\hss}}
\def\Diagonale(#1,#2)#3{\unskip\leavevmode
  \xcoord#1\relax \ycoord#2\relax
      \raise\ycoord \Einheit\hbox to0pt{\hskip\xcoord \Einheit
         \unitlength\Einheit
         \line(1,1){#3}\hss}}
\def\AntiDiagonale(#1,#2)#3{\unskip\leavevmode
  \xcoord#1\relax \ycoord#2\relax \advance\xcoord by -0.05\relax
      \raise\ycoord \Einheit\hbox to0pt{\hskip\xcoord \Einheit
         \unitlength\Einheit
         \line(1,-1){#3}\hss}}
\def\Pfad(#1,#2),#3\endPfad{\unskip\leavevmode
  \xcoord#1 \ycoord#2 \thicklines\ZeichnePfad#3\endPfad\thinlines}
\def\ZeichnePfad#1{\ifx#1\endPfad\let\next\relax
  \else\let\next\ZeichnePfad
    \ifnum#1=1
      \raise\ycoord \Einheit\hbox to0pt{\hskip\xcoord \Einheit
         \vrule height\Pfadd@cke width1 \Einheit depth\Pfadd@cke\hss}%
      \advance\xcoord by 1
    \else\ifnum#1=2
      \raise\ycoord \Einheit\hbox to0pt{\hskip\xcoord \Einheit
        \hbox{\hskip-\PfadD@cke\vrule height1 \Einheit 
width\PfadD@cke depth0pt}\hss}%
      \advance\ycoord by 1
    \else\ifnum#1=3
      \raise\ycoord \Einheit\hbox to0pt{\hskip\xcoord \Einheit
         \unitlength\Einheit
         \line(1,1){1}\hss}
      \advance\xcoord by 1
      \advance\ycoord by 1
    \else\ifnum#1=4
      \raise\ycoord \Einheit\hbox to0pt{\hskip\xcoord \Einheit
         \unitlength\Einheit
         \line(1,-1){1}\hss}
      \advance\xcoord by 1
      \advance\ycoord by -1
    \fi\fi\fi\fi
  \fi\next}
\def\hSSchritt{\leavevmode\raise-.4pt\hbox 
to0pt{\hss.\hss}\hskip.2\Einheit
  \raise-.4pt\hbox to0pt{\hss.\hss}\hskip.2\Einheit
  \raise-.4pt\hbox to0pt{\hss.\hss}\hskip.2\Einheit
  \raise-.4pt\hbox to0pt{\hss.\hss}\hskip.2\Einheit
  \raise-.4pt\hbox to0pt{\hss.\hss}\hskip.2\Einheit}
\def\vSSchritt{\vbox{\baselineskip.2\Einheit\lineskiplimit0pt
\hbox{.}\hbox{.}\hbox{.}\hbox{.}\hbox{.}}}
\def\DSSchritt{\leavevmode\raise-.4pt\hbox to0pt{%
  \hbox to0pt{\hss.\hss}\hskip.2\Einheit
  \raise.2\Einheit\hbox to0pt{\hss.\hss}\hskip.2\Einheit
  \raise.4\Einheit\hbox to0pt{\hss.\hss}\hskip.2\Einheit
  \raise.6\Einheit\hbox to0pt{\hss.\hss}\hskip.2\Einheit
  \raise.8\Einheit\hbox to0pt{\hss.\hss}\hss}}
\def\dSSchritt{\leavevmode\raise-.4pt\hbox to0pt{%
  \hbox to0pt{\hss.\hss}\hskip.2\Einheit
  \raise-.2\Einheit\hbox to0pt{\hss.\hss}\hskip.2\Einheit
  \raise-.4\Einheit\hbox to0pt{\hss.\hss}\hskip.2\Einheit
  \raise-.6\Einheit\hbox to0pt{\hss.\hss}\hskip.2\Einheit
  \raise-.8\Einheit\hbox to0pt{\hss.\hss}\hss}}
\def\SPfad(#1,#2),#3\endSPfad{\unskip\leavevmode
  \xcoord#1 \ycoord#2 \ZeichneSPfad#3\endSPfad}
\def\ZeichneSPfad#1{\ifx#1\endSPfad\let\next\relax
  \else\let\next\ZeichneSPfad
    \ifnum#1=1
      \raise\ycoord \Einheit\hbox to0pt{\hskip\xcoord \Einheit
         \hSSchritt\hss}%
      \advance\xcoord by 1
    \else\ifnum#1=2
      \raise\ycoord \Einheit\hbox to0pt{\hskip\xcoord \Einheit
        \hbox{\hskip-2pt \vSSchritt}\hss}%
      \advance\ycoord by 1
    \else\ifnum#1=3
      \raise\ycoord \Einheit\hbox to0pt{\hskip\xcoord \Einheit
         \DSSchritt\hss}
      \advance\xcoord by 1
      \advance\ycoord by 1
    \else\ifnum#1=4
      \raise\ycoord \Einheit\hbox to0pt{\hskip\xcoord \Einheit
         \dSSchritt\hss}
      \advance\xcoord by 1
      \advance\ycoord by -1
    \fi\fi\fi\fi
  \fi\next}
\def\Koordinatenachsen(#1,#2){\unskip
 \hbox to0pt{\hskip-.5pt\vrule height#2 \Einheit width.5pt depth1 
\Einheit}%
 \hbox to0pt{\hskip-1 \Einheit \xcoord#1 \advance\xcoord by1
    \vrule height0.25pt width\xcoord \Einheit depth0.25pt\hss}}
\def\Koordinatenachsen(#1,#2)(#3,#4){\unskip
 \hbox to0pt{\hskip-.5pt \ycoord-#4 \advance\ycoord by1
    \vrule height#2 \Einheit width.5pt depth\ycoord \Einheit}%
 \hbox to0pt{\hskip-1 \Einheit \hskip#3\Einheit 
    \xcoord#1 \advance\xcoord by1 \advance\xcoord by-#3 
    \vrule height0.25pt width\xcoord \Einheit depth0.25pt\hss}}
\def\Gitter(#1,#2){\unskip \xcoord0 \ycoord0 \leavevmode
  \LOOP\ifnum\ycoord<#2
    \loop\ifnum\xcoord<#1
      \raise\ycoord \Einheit\hbox to0pt{\hskip\xcoord 
\Einheit\Punkt\hss}%
      \advance\xcoord by1
    \repeat
    \xcoord0
    \advance\ycoord by1
  \REPEAT}
\def\Gitter(#1,#2)(#3,#4){\unskip \xcoord#3 \ycoord#4 \leavevmode
  \LOOP\ifnum\ycoord<#2
    \loop\ifnum\xcoord<#1
      \raise\ycoord \Einheit\hbox to0pt{\hskip\xcoord 
\Einheit\Punkt\hss}%
      \advance\xcoord by1
    \repeat
    \xcoord#3
    \advance\ycoord by1
  \REPEAT}
\def\Label#1#2(#3,#4){\unskip \xdim#3 \Einheit \ydim#4 \Einheit
  \def\lo{\advance\xdim by-.5 \Einheit \advance\ydim by.5 \Einheit}%
  \def\llo{\advance\xdim by-.25cm \advance\ydim by.5 \Einheit}%
  \def\loo{\advance\xdim by-.5 \Einheit \advance\ydim by.25cm}%
  \def\o{\advance\ydim by.25cm}%
  \def\ro{\advance\xdim by.5 \Einheit \advance\ydim by.5 \Einheit}%
  \def\rro{\advance\xdim by.25cm \advance\ydim by.5 \Einheit}%
  \def\roo{\advance\xdim by.5 \Einheit \advance\ydim by.25cm}%
  \def\l{\advance\xdim by-.30cm}%
  \def\r{\advance\xdim by.30cm}%
  \def\lu{\advance\xdim by-.5 \Einheit \advance\ydim by-.6 \Einheit}%
  \def\llu{\advance\xdim by-.25cm \advance\ydim by-.6 \Einheit}%
  \def\luu{\advance\xdim by-.5 \Einheit \advance\ydim by-.30cm}%
  \def\u{\advance\ydim by-.30cm}%
  \def\ru{\advance\xdim by.5 \Einheit \advance\ydim by-.6 \Einheit}%
  \def\rru{\advance\xdim by.25cm \advance\ydim by-.6 \Einheit}%
  \def\ruu{\advance\xdim by.5 \Einheit \advance\ydim by-.30cm}%
  #1\raise\ydim\hbox to0pt{\hskip\xdim
     \vbox to0pt{\vss\hbox to0pt{\hss$#2$\hss}\vss}\hss}%
}
\catcode`\@=12

\usepackage{amsmath}
\usepackage{amsthm}
\usepackage[colorlinks=true,
linkcolor=green,
filecolor=brown,
citecolor=green]{hyperref}

\usepackage{color}
\def\red{\textcolor{red} }

\def\s{\ensuremath{\mathcal S}}
\def\si{\ensuremath{\sigma}}
\def\ss{\ensuremath{\s_{[n]}}}

\begin{document}

\begin{center}
{\Large
Counting Stabilized-Interval-Free Permutations         \\ 
}
\vspace{10mm}
DAVID CALLAN  \\
Department of Statistics  \\
University of Wisconsin-Madison  \\
1210 W. Dayton St   \\
Madison, WI \ 53706-1693  \\
{\bf callan@stat.wisc.edu}  \\
\vspace{5mm}
October 10, 2003
\end{center}

\vspace{5mm}

A permutation on $[n]=\{1,2,\ldots,n\}$ is \emph{stabilized-interval-free} 
(SIF) if it does not stabilize any proper subinterval of $[n]$. For 
example, $\left( 
\begin{smallmatrix}
1 & 2 & 3 & 4 & 5 & 6 \\
6 & 1 & 5 & 3 & 4 & 2
\end{smallmatrix} \right) $, or $(3,5,4)(1,6,2)$ in cycle notation, or
$6\,1\,5\,3\,4\,2$ in one-line notation,  fails to be SIF because it 
stabilizes the interval $[3,5]=\{3,4,5\}$. 
On the other hand, the empty permutation is SIF, as is any cycle, and every SIF 
permutation on $[n]$ is fixed-point-free for $n\ge 2$.
Let $a_{n}$ denote the number of SIF permutations on $[n]$ and $A(x)=\sum_{n\ge 
0}^{}a_{n}x^{n}$ their generating function. 
The first objective of this paper is to show that $[x^{n-1}]A(x)^{n}=n!$ 
and hence that the number of SIF 
permutations on $[n]$ is given by 
\htmladdnormallink{A075834}{http://www.research.att.com:80/cgi-bin/access.cgi/as/njas/sequences/eisA.cgi?Anum=A075834}.
This generating function identity amounts to the existence of a 
decomposition of an arbitrary permutation into a list of SIF permutations. The 
second objective is to obtain a recurrence relation that permits 
efficient computation of $a_{n}:$
\[
\ a_{0} =  a_{1} = 1,\ \	
a_{n}  =  \sum_{j=2}^{n-2}(j-1)a_{j}a_{n-j}+(n-1)a_{n-1},\quad n \ge 2.
\]

A little more 
generally, a permutation on a totally ordered set is SIF if it does 
not stabilize any proper saturated chain. Thus  $5\,9\,2\,3$ is SIF 
on $\{2,3,5,9\}$ and its reduced form (replace smallest element by 1, 
second smallest by 2, and so on) is $3\,4\,1\,2$. The former is a 
\emph{labeled} SIF permutation and the latter is \emph{unlabeled}---we take  $[n]$ 
as the standard $n$-element totally ordered set and call a permutation 
on $[n]$ unlabeled; $\s_{[n]}$ denotes the set of all permutations 
on $[n]$. 

For each $\si \in \ss$, one can partition $[n]$ into 
consecutive intervals $I_{1},\ldots,I_{k}$ such that $\si$ 
stabilizes each $I_{j}$. The intervals in the finest such 
partition are called the \emph{components} of $\si$; a permutation with 
exactly one component is \emph{connected} (sometimes called 
indecomposable) 
\htmladdnormallink{A003319}{http://www.research.att.com:80/cgi-bin/access.cgi/as/njas/sequences/eisA.cgi?Anum=A003319}.
Note that the empty 
permutation is not connected.
The restriction of $\si$ to its components clearly gives a 
decomposition of $\si$ into a set of connected permutations on
intervals that partition $[n]$, called the \emph{component permutations} of 
$\si$. These permutations are labeled but we 
also have a decomposition into a \emph{list} of unlabeled
connected permutations of total length $n$ (since we can use position in the 
list to determine the labels) and this decomposition is bijective. For 
example, $3\,2\,5\,1\,4\,7\,8\,6\,9 \longleftrightarrow 
3\,2\,5\,1\,4$\:--\:$2\,3\,1$\:--\:$1$ (the dashes separate list items).

Now $[x^{n-1}]A(x)^{n}$ 
is the number of length-$n$ lists (or simply $n$-lists) of unlabeled SIF permutations of total 
length $n-1$ (keep in mind the empty permutation has length 0). So, to show 
$[x^{n-1}]A(x)^{n}=n!$, it suffices to exhibit a bijection from $\ss$ 
to $n$-lists of unlabeled SIF permutations of total length 
$n-1$, and we will do 
so below. This decomposition into unlabeled SIF permutations is analogous to the 
one above into unlabeled connected permutations but is not so obvious.

Before presenting the bijection we recall some relevant manifestations 
of the Catalan numbers \cite[p.\:219, Ex.\:6.19]{ec2}. A Murasaki diagram is a sequence of vertical 
lines some (all, or none) of which are joined at their tips by 
horizontal lines that never intersect  the interior of a vertical line.
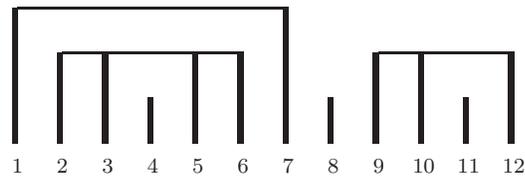
\Einheit=0.6cm
\[
\Pfad(-5,0),222111111\endPfad
\Pfad(-4,0),221111\endPfad 
\Pfad(-3,0),22\endPfad 
\Pfad(-2,0),2\endPfad 
\Pfad(-1,0),22\endPfad
\Pfad(0,0),22\endPfad
\Pfad(1,0),222\endPfad 
\Pfad(2,0),2\endPfad 
\Pfad(3,0),22111\endPfad
\Pfad(4,0),22\endPfad
\Pfad(5,0),2\endPfad 
\Pfad(6,0),22\endPfad 
\Label\u{\scriptstyle{1}}(-5,0)
\Label\u{\scriptstyle{2}}(-4,0)
\Label\u{\scriptstyle{3}}(-3,0)
\Label\u{\scriptstyle{4}}(-2,0)
\Label\u{\scriptstyle{5}}(-1,0)
\Label\u{\scriptstyle{6}}(0,0)
\Label\u{\scriptstyle{7}}(1,0)
\Label\u{\scriptstyle{8}}(2,0)
\Label\u{\scriptstyle{9}}(3,0)
\Label\u{\scriptstyle{10}}(4,0)
\Label\u{\scriptstyle{11}}(5,0)
\Label\u{\scriptstyle{12}}(6,0)
\] 

\noindent The diagram illustrated has 3 components; the first of which has 3 
segments (connected figures), 
the second 1 and the last 2. A partition 
$\{B_{1},B_{2},\ldots,B_{k}\}$ of $[n]$ is noncrossing if $a<b<c<d$ 
with $a,c \in B_{i}$ and $c,d \in B_{j}$ implies $i=j$. Murasaki 
diagrams correspond in an obvious way to noncrossing partitions: the 
one above corresponds to $1\,7$--$2\,3\,5\,6$--$4$--$8$--$9\,10\,12$--$11$ and we 
may speak of the components of a noncrossing partition. 
A lattice path of upsteps $(1,1)$ and downsteps $(1,-1)$ (starting at 
the origin for convenience) is balanced if it ends on the $x$-axis, 
nonnegative if it never dips below the $x$-axis, Dyck if it is both. A 
Dyck $n$-path $P$ has $n$ upsteps and $n$ downsteps; each downstep $d$ 
has a matching upstep $u$: head  horizontally west from $d$ to the 
first upstep $u$ that you encounter.  Each $x$-axis 
point on $P$ other than the starting point is a return of $P$; $P$
is strict if it has only one return. Its returns divide a nonempty Dyck path 
into a list of its components, each of which is a strict Dyck path. 
For any path, a nonzero ascent is a 
maximal sequence of contiguous upsteps (we assume a zero ascent 
between a pair of contiguous downsteps); similarly for descents.

Noncrossing partitions $\pi$ on $[n]$ correspond to Dyck $n$-paths $P$: 
arrange the blocks of $\pi$ in increasing order of their maximal 
elements; let $(m_{i})_{i=1}^{k}$ be these maximal elements and let 
$(n_{i})_{i=1}^{k}$ be the corresponding block sizes. Then, with 
$m_{0}:=0$, the lists $(m_{i}-m_{i-1})_{i=1}^{k}$  and 
$(n_{i})_{i=1}^{k}$ determine $\pi$ and are, respectively, the nonzero 
ascent lengths and nonzero descent lengths defining $P$. This 
correspondence preserves components.

An arbitrary permutation $\si$ can be split into a set of 
labeled SIF permutations whose underlying sets partition $[n]$.
First, decompose $\si$ into its connected components 
$\si_{1},\si_{2},\ldots,\si_{k}$. Set aside all stabilized proper
subintervals (if any) of each $\si_{i}$; what's left will be $k$ 
nonempty SIF permutations. Repeat this procedure on the entire 
permutation that was set aside, continuing till nothing is set aside.
The resulting set of SIFs corresponds to a Murasaki diagram 
in which an unlabeled SIF is associated with each segment; the 
segments record the underlying sets, the unlabeled SIFs the action of 
the permutation.
For example, 
$\left(
\begin{smallmatrix}
1 & 2 & 3 & 4 & 5 & 6 & 7 & 8 & 9 & 10 & 11 & 12 \\
7 & 5 & 6 & 4 & 2 & 3 & 1 & 8 & 10 & 12 & 11 & 9
\end{smallmatrix} \right)$ 
splits into 
$\left(
\begin{smallmatrix}
1  & 7  \\
7 & 1 
\end{smallmatrix} \right),\, \left(
\begin{smallmatrix}
 2 & 3 & 5 & 6  \\
 5 & 6 & 2 & 3 
\end{smallmatrix} \right),\, \left(
\begin{smallmatrix}
 4  \\
 4 
\end{smallmatrix} \right),\, \left(
\begin{smallmatrix}
 8  \\
 8 
\end{smallmatrix} \right),\, \left(
\begin{smallmatrix}
 9 & 10  & 12 \\
10 & 12  & 9
\end{smallmatrix} \right),\, \left(
\begin{smallmatrix}
 11  \\
 11 
\end{smallmatrix} \right)$. 
The Murasaki diagram is the one above and 
unlabeled SIFs are associated with segments as follows.
\[
\begin{array}{ccccccc}
\textrm{segments by smallest element} & 1 & 2 & 4 & 8 & 9 & 11  \\ 
\hline
\textrm{corresponding unlabeled SIF } & 21 & 3412 & 1 & 1 & 231 & 1
\end{array}
\]

Now we are ready to present the bijection from $\ss$ to  
$n$-lists of unlabeled SIF permutations whose total length is
$n-1$, and we will use
\[
\setcounter{MaxMatrixCols}{16}
\si=\left(
\begin{matrix}
	1 & 2 & 3 & 4 & 5 & 6 & 7 & 8 & 9 & 10 & 11 & 12 & 13 & 14 & 15 & 16  \\
	2 & 4 & 3 & 1 & 8 & 7 & 6 & 5 & 13 & 10 & 9 & 16 & 11 & 15 & 14 & 12
\end{matrix}
\right)
\]
as a working example with $n=16$. First, decompose $\si$ into its 
components $(\si_{i})_{i=1}^{k}$; note that $n$ will occur in the 
last one. Record the 
position $j$ of $n$ in $\si$ (here $j=12$), then delete $n$ from $\si_{k}$ to get a 
permutation $\si_{k}'$ 
(deleting $n$ simply means erasing $n$ from its cycle and so 
$\si_{k}'(j)=\si(n)$\:) and $j$ is necessarily in the first component 
of $\si_{k}'$ because $\si_{k}$ is connected. Now draw the Murasaki diagrams for 
$\si_{1},\ldots,\si_{k-1},\si_{k}'$ and 
record the associated unlabeled SIF for each segment. 
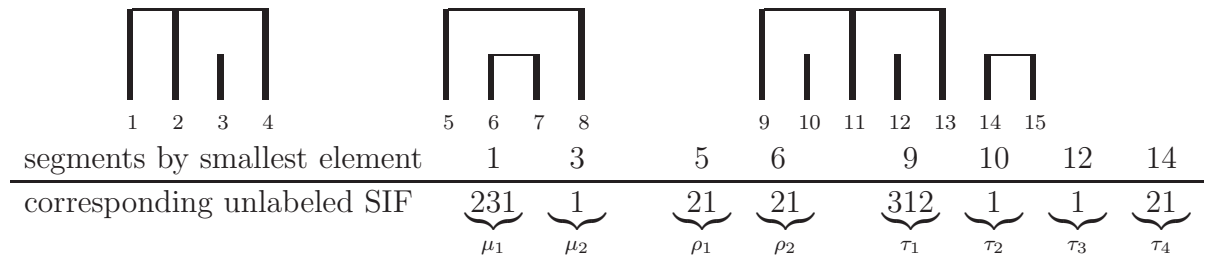
\Einheit=0.6cm
\[
\Pfad(-10,0),22111\endPfad
\Pfad(-9,0),22\endPfad 
\Pfad(-8,0),2\endPfad 
\Pfad(-7,0),22\endPfad 
\Pfad(-3,0),22111\endPfad
\Pfad(-2,0),21\endPfad 
\Pfad(-1,0),2\endPfad 
\Pfad(-0,0),22\endPfad 
\Pfad(4,0),221111\endPfad
\Pfad(5,0),2\endPfad 
\Pfad(6,0),22\endPfad 
\Pfad(7,0),2\endPfad 
\Pfad(8,0),22\endPfad 
\Pfad(9,0),21\endPfad
\Pfad(10,0),2\endPfad 
\Label\u{\scriptstyle{1}}(-10,0)
\Label\u{\scriptstyle{2}}(-9,0)
\Label\u{\scriptstyle{3}}(-8,0)
\Label\u{\scriptstyle{4}}(-7,0)
\Label\u{\scriptstyle{5}}(-3,0)
\Label\u{\scriptstyle{6}}(-2,0)
\Label\u{\scriptstyle{7}}(-1,0)
\Label\u{\scriptstyle{8}}(-0,0)
\Label\u{\scriptstyle{9}}(4,0)
\Label\u{\scriptstyle{10}}(5,0)
\Label\u{\scriptstyle{11}}(6,0)
\Label\u{\scriptstyle{12}}(7,0)
\Label\u{\scriptstyle{13}}(8,0)
\Label\u{\scriptstyle{14}}(9,0)
\Label\u{\scriptstyle{15}}(10,0)
\] 
\[
\begin{array}{ccccccccccc}
\textrm{segments by smallest element} \hspace*{2mm} & 1 & 3 & 
\hspace*{2mm} & 5 & 
6 \hspace*{2mm} & & 9 & 10 & 12 & 14  \\ 
\hline
\textrm{corresponding unlabeled SIF } \hspace*{2mm} & 
\underbrace{231}_{\mu_{1}} & \underbrace{1}_{\mu_{2}} &  
\hspace*{2mm} &
\underbrace{21}_{\rho_{1}} & \underbrace{21}_{\rho_{2}} & \hspace*{2mm} & 
\underbrace{312}_{\tau_{1} } & \underbrace{1}_{ \tau_{2}}  & \underbrace{1}_{\tau_{3}}  & 
\underbrace{21}_{\tau_{4}} 
\end{array}
\]

Translate each Murasaki diagram $\rightarrow$ noncrossing partition $\rightarrow$ 
Dyck path, recalling that segment $\rightarrow$ block $\rightarrow$ 
nonzero descent, so each nonzero descent is associated with an SIF, and mark upstep $j$ 
(in red in the following figure) unless $j=n$ in which case 
$\si_{k}'$ is the empty permutation.
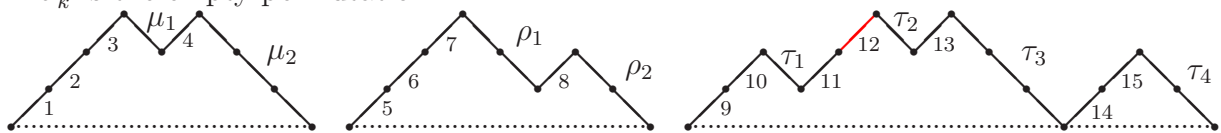
\Einheit=0.5cm
\[
\Pfad(-16,0),33343444\endPfad
\Pfad(-7,0),33344344\endPfad
\Pfad(2,0),3343\endPfad
\red{\Pfad(6,2),3\endPfad}
\Pfad(7,3),434443344\endPfad
\SPfad(-16,0),11111111\endSPfad
\SPfad(-7,0),11111111\endSPfad
\SPfad(2,0),11111111111111\endSPfad
\DuennPunkt(-16,0)
\DuennPunkt(-15,1)
\DuennPunkt(-14,2)
\DuennPunkt(-13,3)
\DuennPunkt(-12,2)
\DuennPunkt(-11,3)
\DuennPunkt(-10,2)
\DuennPunkt(-9,1)
\DuennPunkt(-8,0)
\DuennPunkt(-7,0)
\DuennPunkt(-6,1)
\DuennPunkt(-5,2)
\DuennPunkt(-4,3)
\DuennPunkt(-3,2)
\DuennPunkt(-2,1)
\DuennPunkt(-1,2)
\DuennPunkt(0,1)
\DuennPunkt(1,0)
\DuennPunkt(2,0)
\DuennPunkt(3,1)
\DuennPunkt(4,2)
\DuennPunkt(5,1)
\DuennPunkt(6,2)
\DuennPunkt(7,3)
\DuennPunkt(8,2)
\DuennPunkt(9,3)
\DuennPunkt(10,2)
\DuennPunkt(11,1)
\DuennPunkt(12,0)
\DuennPunkt(13,1)
\DuennPunkt(14,2)
\DuennPunkt(15,1)
\DuennPunkt(16,0)
\Label\o{\scriptstyle{1}}(-15.0,0)
\Label\o{\scriptstyle{2}}(-14.3,0.7)
\Label\o{\scriptstyle{3}}(-13.3,1.7)
\Label\o{\scriptstyle{4}}(-11.3,1.8)
\Label\o{\scriptstyle{5}}(-6.0,0.0)
\Label\o{\scriptstyle{6}}(-5.3,0.7)
\Label\o{\scriptstyle{7}}(-4.3,1.7)
\Label\o{\scriptstyle{8}}(-1.3,0.7)
\Label\o{\scriptstyle{9}}(3.0,0.0)
\Label\o{\scriptstyle{10}}(3.8,0.7)
\Label\o{\scriptstyle{11}}(5.8,0.7)
\Label\o{\scriptstyle{12}}(6.8,1.7)
\Label\o{\scriptstyle{13}}(8.8,1.7)
\Label\o{\scriptstyle{14}}(13.0,-0.1)
\Label\o{\scriptstyle{15}}(13.8,0.7)
\Label\o{\mu_{1}}(-12.0,2.3)
\Label\o{\mu_{2}}(-8.8,1.4)
\Label\o{\rho_{1}}(-2.2,1.9)
\Label\o{\rho_{2}}(.7,1.0)
\Label\o{\tau_{1}}(4.8,1.3)
\Label\o{\tau_{2}}(7.8,2.3)
\Label\o{\tau_{3}}(11.2,1.4)
\Label\o{\tau_{4}}(15.6,.9)
\]
\begin{center}
\small{Dyck paths, upsteps labeled in order, nonzero descents labeled with 
corresponding SIF permutation (its length = length of descent). The matching 
upsteps for a nonzero descent give a block of the noncrossing partition and 
identify a segment of the Murasaki diagram. All but the last 
are strict Dyck paths. The marked upstep (here 12) is in the first component of 
the last path (unless the last path is empty).}
\end{center}

 We can use a cut-and-paste technique to massage these Dyck 
paths into a balanced path in a reversible way (making critical use 
of the marked upstep). The process will preserve all nonzero descents 
and so we can carry their SIF labels along with them. Cut the last 
Dyck path just before its marked upstep into two paths $R,S$.
For each preceding Dyck 
path, remove its last upstep thereby forming a path $P_{i}$ and a 
nonzero descent $D_{i},\ 1\le i \le k-1$, and $k-1$ removed upsteps 
$u$. Then rearrange in the following order to form a balanced path 
$Q$: $D_{1}\ u\ D_{2}\ u\ \ldots 
\ D_{k-1}\ u\ S\ R \ P_{1}\ P_{2}\ \ldots \ P_{k-1}$.
\Einheit=0.5cm
\[
 \Pfad(-15,5),444344334344433443343333433344\endPfad
\SPfad(-15,5),111111111111111111111111\endSPfad
\SPfad(11,5),1111\endSPfad
\DuennPunkt(-15,5)
\DuennPunkt(-15,5)
\DuennPunkt(-14,4)
\DuennPunkt(-13,3)
\DuennPunkt(-12,2)
\DuennPunkt(-11,3)
\DuennPunkt(-10,2)
\DuennPunkt(-9,1)
\DuennPunkt(-8,2)
\DuennPunkt(-7,3)
\DuennPunkt(-6,2)
\DuennPunkt(-5,3)
\DuennPunkt(-4,2)
\DuennPunkt(-3,1)
\DuennPunkt(-2,0)
\DuennPunkt(-1,1)
\DuennPunkt(0,2)
\DuennPunkt(1,1)
\DuennPunkt(2,0)
\DuennPunkt(3,1)
\DuennPunkt(4,2)
\DuennPunkt(5,1)
\DuennPunkt(6,2)
\DuennPunkt(7,3)
\DuennPunkt(8,4)
\DuennPunkt(9,5)
\DuennPunkt(10,4)
\DuennPunkt(11,5)
\DuennPunkt(12,6)
\DuennPunkt(13,7)
\DuennPunkt(14,6)
\DuennPunkt(15,5)
\Label\o{\scriptstyle{4}}(-11.3,1.8)
\Label\o{\scriptstyle{8}}(-8.2,0.8)
\Label\o{\scriptstyle{12}}(-7.1,1.8)
\Label\o{\scriptstyle{13}}(-5.1,1.8)
\Label\o{\scriptstyle{14}}(-1.1,-0.2)
\Label\o{\scriptstyle{15}}(-0.1,0.8)
\Label\o{\scriptstyle{9}}(2.9,-0.2)
\Label\o{\scriptstyle{10}}(3.9,0.7)
\Label\o{\scriptstyle{11}}(5.9,0.7)
\Label\o{\scriptstyle{1}}(6.7,1.7)
\Label\o{\scriptstyle{2}}(7.7,2.7)
\Label\o{\scriptstyle{3}}(8.7,3.7)
\Label\o{\scriptstyle{5}}(10.7,3.7)
\Label\o{\scriptstyle{6}}(11.8,4.9)
\Label\o{\scriptstyle{7}}(12.7,5.7)
\Label\o{\mu_{2}}(-12.9,3.3)
\Label\o{\rho_{2}}(-9.4,1.9)
\Label\o{\tau_{2}}(-6.2,2.3)
\Label\o{\tau_{3}}(-2.9,1.4)
\Label\o{\tau_{4}}(1.3,.9)
\Label\o{\tau_{1}}(4.8,1.3)
\Label\o{\mu_{1}}(10.0,4.3)
\Label\o{\rho_{1}}(14.4,6.0)
\Label\u{Q}(0,-.6)
\]
\begin{center}
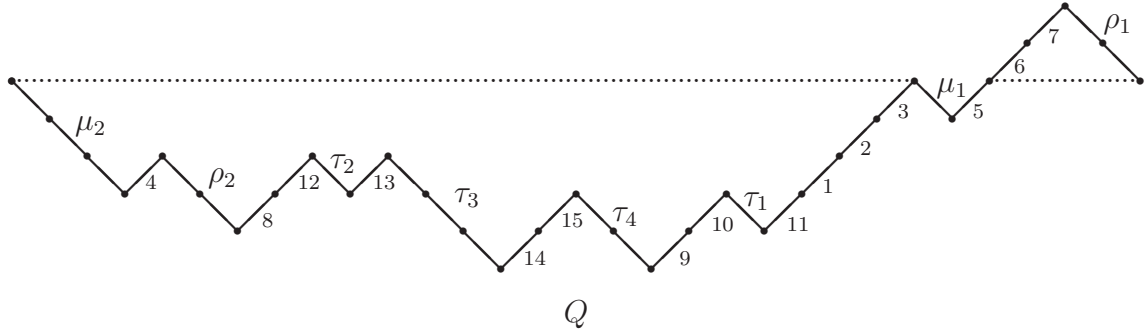
\small{balanced path, descents labeled with 
corresponding SIF permutation,\\
upsteps shown with their original numbering 
to aid the reader}
\end{center}  

The original Dyck paths can be recovered from the balanced 
path. In brief, the center point $p$ of the first double rise (= 
consecutive pair of upsteps) identifies the initial 
vertex of the marked upstep. The path from $p$ to the 
rightmost lowest point of $Q$ following $p$ is $S$ (this relies on the fact 
that the marked upstep 
was in the first component); from there to the rightmost point 
$q$ at $p$'s level is $R$. The descents preceding $p$ are 
$D_{1},\ldots,D_{k-1}$ and their lengths determine how 
far to proceed from $q$ to recover $P_{1},\ldots,P_{k-1}$.

The preceding outline needs a little elaboration to cover special cases. More 
precisely, prepend and append upsteps to $Q$ to guarantee the 
existence of a double rise and the point $p$. If $Q$ starts with an 
upstep, then $p$ will be the origin, the list $D_{1},\ldots,D_{k-1}$ 
will be vacuous and the original permutation $\si$ 
will be connected. If $p$ is the last point of $Q$, then $Q$ will have 
a sawtooth shape, $\backslash/\backslash/\backslash/\backslash/$, 
and $\si=$ identity. If $n$ is a fixed point of $\si$, then the last 
Dyck path is empty (there is no marked upstep) and $Q$ proceeds from 
$p$ with an upstep and never drops back to the level of $p$. Also, of 
course, either one of the paths $R,S$ may be empty.

Finally, scan all descents of the balanced path $Q$, recording 
$\emptyset$ (the empty permutation) for each zero descent and its associated 
unlabeled SIF permutation for each nonzero descent.
\[\underbrace{
\ \mu_{2}\ \:\rho_{2}\ \:\emptyset\ \:\tau_{2}\ \:\tau_{3}\
\emptyset\ \:\tau_{4}\ \:\emptyset\ \:\tau_{1}\ \:\emptyset\ 
\:\emptyset\ \:\emptyset\
\mu_{1}\ \:\emptyset\ \:\emptyset\ \:\rho_{1}\ }_{\textrm{\normalsize{ 
$n$-list of SIF
permutations of total length $n-1$}}}
\]

We have shown that the generating function for the number $a_{n}$ of 
SIF permutations on $[n]$ is that of 
\htmladdnormallink{A075834}{http://www.research.att.com:80/cgi-bin/access.cgi/as/njas/sequences/eisA.cgi?Anum=A075834}
but to calculate values of $a_{n}$ it is more efficient to develop a 
recurrence relation. Let $a_{n,k}$ denote the number of
permutations on $[n]$ that do not stabilize any proper subinterval beginning 
at $i$ for $i<k$. Thus $a_{n,1} = n!$
\htmladdnormallink{A000142}{http://www.research.att.com:80/cgi-bin/access.cgi/as/njas/sequences/eisA.cgi?Anum=A000142}, 
$a_{n,2}$ is the number of connected permutations on $[n]$ 
\htmladdnormallink{A003319 }{http://www.research.att.com:80/cgi-bin/access.cgi/as/njas/sequences/eisA.cgi?Anum=A003319 }
(apart from the first term---we need to set $a_{1,2}=0$), and $a_{n,n}=a_{n}$. 
Counting permutations by their first  stabilized subinterval, it is 
straightforward to obtain 
the following recurrence (given in Mathematica code ).
\begin{verbatim}
c[0]=0; c[n_]/;n>=1 := c[n] = n!-Sum[c[i](n-i)!,{i,n-1}] 
(* c[n] = # connected perms on [n] *)
a[n_,k_]/;n>=0  && k==n+1 := 0; 
a[n_,1]/;n>=1 := n!; 
a[n_,k_]/;2<=k<=n := a[n,k] = 
   n!-Sum[c[j-i+1]a[n-(j-i+1),i],{i,k-1},{j,i,n}];	
\end{verbatim}

However, there is also a direct recurrence for $a_{n}$ (vacuous sums 
are 0):
\[
\ a_{0} =  a_{1} = 1,\ \	
a_{n}  =  \sum_{j=2}^{n-2}(j-1)a_{j}a_{n-j}+(n-1)a_{n-1},\quad n \ge 2.
\]
The right hand side above counts SIF permutations $\si$ on $[n]$ by 
the parameter $j=n-1-s$ where $s$ is the size of the largest proper 
subinterval $I$ of $[n-1]$ such that $\si$ stabilizes $I \cup 
\{n\}$. ($I$ is necessarily an interior interval of $[n-1]$ and may be 
empty.)

To see this, first note that if $\si_{n-1}$ is SIF on $[n-1]$ and $n$ 
is inserted anywhere into a cycle of $\si_{n-1}\ (n-1$ possible ways) to 
form $\si \in \ss$, then $\si$ is also SIF. This accounts for the last 
term. Now suppose $\si $ is SIF on $[n]$ and the result $\si_{n-1}\in 
\mathcal{S}_{[n-1]}$ of deleting $n$ from its cycle in $\si$ fails to 
be SIF. Consider the maximal proper subintervals of $[n-1]$ 
stabilized by $\si_{n-1}$. There is at least one such by assumption and 
at most one, call it $I$, because otherwise $\si$ itself would 
stabilize all but one of them. Let $\rho$ denote the 
restriction of $\si_{n-1}$ to $I$ and $\tau$ the restriction of 
$\si_{n-1}$ to $[n-1]\backslash I$. Then $\si$ is obtained from the 
pair $\tau,\rho$ by inserting $n$ into a cycle of $\rho$, not $\tau$, 
otherwise $\si$ would stabilize $I$. We may write the interval $I$ as 
$[k+1,n-j+k-1]$ for some $1 \le k < j \le n-2$ so that the size of $I$ 
is $s:=n-j-1$ and $I$ is clearly the largest proper subinterval of 
$[n-1]$ such that $\si$ stabilizes $I \cup 
\{n\}$. Now $\tau$ is SIF on $[n-1]\backslash I$ by definition of 
$\rho$. We claim $\rho':=\si$  restricted to $I \cup \{n\}$ is SIF 
also: if $\rho'$ stabilized a proper subinterval of $I$, then $\si$ 
would too, and if $\rho'$ stabilized a proper terminal subinterval 
(containing $n$), then $\si$ would stabilize the corresponding initial 
subinterval. All told, for each $j \in [2,n-2]$, we have $j-1$ choices 
for $k$ and every permutation $\si$ formed in this way from SIF 
permutations $\rho'$ on $I \cup \{n\}$ ($a_{n-j}$ choices) and $\tau$ 
on $[n-1]\backslash I$ $(a_{j}$ choices) is SIF. The recurrence 
follows. We note that it implies the differential equation
\[
xA'(x)=A(x)-x-\frac{x}{A(x)-1}
\]
for the generating function $A(x)$.

Asymptotically, the proportion of permutations on $[n]$ that are connected
(indecomposable) is $1-\frac{2}{n}+O(\frac{1}{n^{2}})$ \cite[p.\:295,
Ex.\:16]{comtet} and there is a simple heuristic
explanation: the easiest way for a permutation on $[n]$ to be
decomposable is for it to fix $1$ or $n$ and there are
$2(n-1)!-(n-2)!$ permutations that do so. Far fewer permutations
stabilize any other initial interval and so the dominant term in the
number of decomposable permutations on $[n]$ is $2(n-1)!$. Similarly,
the easiest way for $\si \in \ss$ to fail to be SIF is for it to have
a fixed point. The proportion of fixed-point-free permutations on $[n]$
is well known to be very near $\frac{1}{e}$, suggesting that the proportion
of SIF permutations on $[n]$ is $\frac{1}{e}
+O(\frac{1}{n})$, and indeed computer calculations suggest it is
$\frac{1}{e}(1-\frac{1}{n})+O(\frac{1}{n^{2}})$ and maybe 
$\frac{1}{e}(1-\frac{1}{n}-\frac{5}{2n^{2}})+O(\frac{1}{n^{3}})$. It would be
interesting to prove this.

\end{document}